\documentclass{conm-p-l}

\usepackage{amsmath,amsthm,amsfonts,amscd,amssymb}

\newtheorem{thm}{Theorem}[section]

\newtheorem{conj}[thm]{Conjecture}

\newtheorem{prob}{Problem}

\theoremstyle{definition}

\theoremstyle{remark}
\newtheorem{rem}{Remark}[section]

\begin{document}

\title[Heights and quadratic forms]{Heights and quadratic forms: Cassels' theorem and its generalizations}
\author{Lenny Fukshansky}

\address{Department of Mathematics, 850 Columbia Avenue, Claremont McKenna College, Claremont, CA 91711}
\email{lenny@cmc.edu}

\begin{abstract}
In this survey paper, we discuss the classical Cassels' theorem on existence of small-height zeros of quadratic forms over $\mathbb Q$ and its many extensions, to different fields and rings, as well as to more general situations, such as existence of totally isotropic small-height subspaces. We also discuss related recent results on effective structural theorems for quadratic spaces, as well as Cassels'-type theorems for small-height zeros of quadratic forms with additional conditions. We conclude with a selection of open problems.
\end{abstract}

\maketitle

\def\A{{\mathcal A}}
\def\AA{{\mathfrak A}}
\def\B{{\mathcal B}}
\def\C{{\mathcal C}}
\def\D{{\mathcal D}}
\def\F{{\mathcal F}}
\def\x{{\mathcal H}}
\def\I{{\mathcal I}}
\def\J{{\mathcal J}}
\def\K{{\mathcal K}}
\def\kk{{\mathfrak K}}
\def\L{{\mathcal L}}
\def\M{{\mathcal M}}
\def\W{{\omega}}
\def\mm{{\mathfrak m}}
\def\MM{{\mathfrak M}}
\def\NN{{\mathfrak N}}
\def\OO{{\mathfrak O}}
\def\O{{\mathcal O}}
\def\R{{\mathcal R}}
\def\s{{\mathcal S}}
\def\V{{\mathcal V}}
\def\UU{{\mathfrak U}}
\def\X{{\mathcal X}}
\def\Y{{\mathcal Y}}
\def\Z{{\mathcal Z}}
\def\H{{\mathcal H}}
\def\G{{\mathcal G}}
\def\cee{{\mathbb C}}
\def\pee{{\mathbb P}}
\def\que{{\mathbb Q}}
\def\real{{\mathbb R}}
\def\zed{{\mathbb Z}}
\def\aaa{{\mathbb A}}
\def\hyp{{\mathbb H}}
\def\Bb{{\mathbb B}}
\def\ff{{\mathbb F}}
\def\Nn{{\mathbb N}}
\def\kk{{\mathfrak K}}
\def\qbar{{\overline{\mathbb Q}}}
\def\kbar{{\overline{K}}}
\def\xbar{{\overline{x}}}
\def\ybar{{\overline{Y}}}
\def\kkbar{{\overline{\mathfrak K}}}
\def\ubar{{\overline{U}}}
\def\abar{{\overline{a}}}
\def\eps{{\varepsilon}}
\def\ahat{{\hat \alpha}}
\def\bhat{{\hat \beta}}
\def\gt{{\tilde \gamma}}
\def\h{{\tfrac12}}
\def\ba{{\boldsymbol a}}
\def\be{{\boldsymbol e}}
\def\bei{{\boldsymbol e_i}}
\def\bc{{\boldsymbol c}}
\def\bm{{\boldsymbol m}}
\def\bk{{\boldsymbol k}}
\def\bi{{\boldsymbol i}}
\def\bl{{\boldsymbol l}}
\def\bq{{\boldsymbol q}}
\def\bu{{\boldsymbol u}}
\def\bt{{\boldsymbol t}}
\def\bs{{\boldsymbol s}}
\def\bv{{\boldsymbol v}}
\def\bw{{\boldsymbol w}}
\def\bx{{\boldsymbol x}}
\def\bX{{\boldsymbol X}}
\def\bz{{\boldsymbol z}}
\def\bwy{{\boldsymbol y}}
\def\bY{{\boldsymbol Y}}
\def\bL{{\boldsymbol L}}
\def\bet{{\boldsymbol\eta}}
\def\bxi{{\boldsymbol\xi}}
\def\bo{{\boldsymbol 0}}
\def\bol{{\boldkey 1}_L}
\def\ep{\varepsilon}
\def\p{\boldsymbol\varphi}
\def\q{\boldsymbol\psi}
\def\Hf{H_{\fin}^{\O}}
\def\Hfo{H_{\fin}^{\O_1}}
\def\Hft{H_{\fin}^{\O_2}}
\def\Hfd{H_{\fin}^{O_D}}
\def\rank{\operatorname{rank}}
\def\aut{\operatorname{Aut}}
\def\lcm{\operatorname{lcm}}
\def\sgn{\operatorname{sgn}}
\def\spn{\operatorname{span}}
\def\md{\operatorname{mod}}
\def\Norm{\operatorname{Norm}}
\def\dim{\operatorname{dim}}
\def\det{\operatorname{det}}
\def\Vol{\operatorname{Vol}}
\def\rk{\operatorname{rk}}
\def\ord{\operatorname{ord}}
\def\ker{\operatorname{ker}}
\def\div{\operatorname{div}}
\def\Gal{\operatorname{Gal}}
\def\Tr{\operatorname{Tr}}
\def\nn{\operatorname{N}}
\def\inf{\operatorname{inf}}
\def\fin{\operatorname{fin}}
\def\Gr{\operatorname{Gr}}
\def\Mat{\operatorname{Mat}}
\def\GL{\operatorname{GL}}
\def\disc{\operatorname{disc}}

\section{Introduction: Cassels' theorem}
\label{intro}

Given a symmetric bilinear form 
\begin{equation}
\label{form}
F(\bX,\bY) = F(X_1,\dots,X_N,Y_1,\dots,Y_N) := \sum_{i=1}^N \sum_{j=1}^N f_{ij} X_i Y_j
\end{equation}
in $2N \geq 4$ variables with rational coefficients $f_{ij}=f_{ji}$, it is a classical problem to determine whether the corresponding quadratic form $F(\bX) := F(\bX,\bX)$ in $N \geq 2$ variables is isotropic over $\que$. The answer is famously provided by the Hasse-Minkowski theorem: $F$ is isotropic over $\que$ if and only if it is isotropic over every completion of $\que$. A corollary of this result is Meyer's theorem, which guarantees that if $F$ is indefinite and $N \geq 5$, then $F$ is necessarily isotropic over $\que$. Both of these results however are {\it ineffective} in the sense that they do not provide any insight into how to find a nontrivial zero of $F$ over $\que$, should one exist. One possible {\it effective} approach to this problem would be the following. Suppose that we can prove that whenever $F$ is isotropic over $\que$, there must exist a nonzero vector $\bx \in \zed^N$ such that $F(\bx) = 0$ and
\begin{equation}
\label{sb1}
\max_{1 \leq i \leq N} |x_i| \leq C_N(F),
\end{equation}
where $C_N(F)$ is some explicit expression depending on $F$ and $N$. Since the set of points
$$\left\{ \bx \in \zed^N : \max_{1 \leq i \leq N} |x_i| \leq C_N(F) \right\}$$
is finite, one can now look for a nontrivial integral zero of $F$ by searching through this set. This consideration motivates calling the explicit expression $C_N(F)$ in the upper bound of \eqref{sb1} a {\it search bound} for $F$ over $\que$. Homogeneity of $F(\bX)$ implies that we can assume without loss of generality that it has integer coefficients and look for integral zeros. In his celebrated 1955 paper \cite{cassels:small}, J. W. S. Cassels found such a search bound.

\begin{thm}[Cassels' Theorem] \label{cassels}Let 
$$F(\bX) = \sum_{i=1}^N \sum_{j=1}^N f_{ij} X_i X_j \in \zed[\bX]$$
be an isotropic integral quadratic form in $N \geq 2$ variables, then there exists $\bx \in \zed^N \setminus \{\bo\}$ such that $F(\bx)=0$ and 
\begin{equation}
\label{cassels_bound}
\max_{1 \leq i \leq N} |x_i|  \leq \left( 3 \sum_{i=1}^N \sum_{j=1}^N |f_{ij}| \right)^{\frac{N-1}{2}}.
\end{equation}
\end{thm}

Our presentation of Cassels' theorem follows not the original version of \cite{cassels:small}, but rather a refined version recorded in Section~6.8 of \cite{cassels_book}. In fact, the exponent in the upper bound of \eqref{cassels_bound} is best possible in general, as demonstrated by an example due to M. Kneser (see p.~87 of \cite{cassels_book}). On the other hand, G. L. Watson showed \cite{watson} that the exponent can be replaced by $\max \{ 2,r/2,s/2 \}$, where $r$ is the number of positive squares and $s$ the number of negative squares when $F(\bX)$ is expressed as a sum of squares (with signs) of real linear forms. More recently, T. D. Browning and R. Dietmann~\cite{browning_dietmann} showed that it is possible to obtain smaller exponents in the upper bound~\eqref{cassels_bound} for {\it generic} quadratic forms.

Cassels' original argument proceeds as follows. Let $\bx$ be the smallest (with respect to sup-norm)  nontrivial integral zero of $F$. Minkowski's Linear Forms Theorem is then applied to construct an anisotropic integral point $\bwy$ of small sup-norm. If $\bx$ does not satisfy \eqref{cassels_bound}, then its orthogonal reflection with respect to $F$ in the hyperplane orthogonal to $\bwy$ has smaller sup-norm then $\bx$, leading to a contradiction. A different proof of the same result (although with a different constant in the upper bound) was exhibited by H. Davenport \cite{davenport1} and then generalized by B. Birch and H. Davenport \cite{birch_davenport}. 

Cassels' theorem can be viewed as one of the rare instances of effective search bounds for Diophantine equations. Indeed, suppose that we have a general polynomial $F$ of degree $M$ in $N$ variables with integer coefficients, then Hilbert's tenth problem asks whether there exists an algorithm to determine if $F$ has a nontrivial integral zero. Assuming we knew that whenever such a zero exists, there must exist one satisfying \eqref{sb1} with an explicit bound $C_N(F)$, Hilbert's tenth problem could be reduced to a finite search algorithm. Therefore Matiyasevich's famous negative answer \cite{mat} to Hilbert's tenth problem indicates that search bounds in general cannot exist. Moreover, it has been shown by J. P. Jones \cite{jp_jones} that the question whether a general Diophantine equation of degree four or larger has a solution in positive integers is already undecidable, and not much else is known for polynomials of degree $\geq 4$. Although some partial results are available for cubic polynomials (see \cite{bde} for the most recent developments), the general question of obtaining search bounds for solutions of an arbitrary integral cubic polynomial is open. Hence the only two well understood cases for existence of search bounds are those of linear and quadratic polynomials. The linear case is resolved by Siegel's lemma and its numerous generalizations (see \cite{null} for a recent account), while the quadratic case is addressed by Cassels' theorem and a variety of results extending it. It should also be remarked that Cassels' theorem and related results have been successfully applied to the problems of solubility of Diophantine inequalities~\cite{birch_davenport} and to bounding size of representations of numbers by quadratic forms~\cite{chalk}.

In this note we present a survey of extensions of Cassels' theorem with a view towards a unified presentation and generality. All these extensions and generalizations can be split into several categories:

\begin{enumerate}

\item Extensions to quadratic forms over more general coefficient fields and rings, such as number fields, function fields, and quaternion algebras.

\item Extensions to statements about multiple linearly independent zeros of a quadratic form and totally isotropic subspaces of quadratic spaces.

\item Effective structural results for quadratic spaces over global fields, such as effective versions of Witt decomposition theorem and Cartan-Dieudonn\'e theorem on decomposition of isometries.

\item Effective results on zeros of quadratic forms satisfying additional algebraic conditions and analogues of Cassels' theorem for inhomogeneous quadratic equations.

\end{enumerate}

\noindent
The results we discuss here are stated in terms of {\it height functions}, which are natural global-field analogues of the sup-norm on integers from the stand-point of search bounds. In other words,  \eqref{cassels_bound} can be viewed as a bound on the height of a nontrivial zero $\bx$ of the quadratic form $F$ in terms of the height of this form. All of the bounds we review are effective, however in the interest of simplicity of presentation we do not write out the constants explicitly, just their dependences (explicit values of the constants can be found in the references we provide). The important feature of the bounds which we try to emphasize here is their dependence on the heights of quadratic forms and spaces with explicit exponents, which are polynomial in the dimension of the space in a majority of situations. This paper is structured as follows. In Section~\ref{heights} we introduce the necessary notation and the machinery of height functions. We then review progress in each of the four categories outlined above in Sections~\ref{fields}--\ref{inhomog}. We conclude with a selection of general open problems and directions in Section~\ref{problems}.

\begin{rem} \label{survey} A survey of progress on the general problem of finding bounds for zeros of a rational quadratic form in terms of its coefficients has been published in 1990 by H. P. Schlickewei and W. M. Schmidt~\cite{schm_schlick}. There have however been substantial further developments in the subject since this paper appeared. It is the goal of this note to review the problem from its beginnings, especially concentrating on the recent work (since Schlickewei - Schmidt review) and presenting a unified approach of the effective theory of quadratic forms via height functions.
\end{rem} 
\bigskip

\section{Notation and heights}
\label{heights}

We start with some notation.  Our presentation of absolute values and heights is unified over several different types of fields, largely following~\cite{absolute:siegel}. Throughout this paper, $K$ will either be a number field (a finite extension of $\que$) or a function field (a finite extension of the field $\kk = \kk_0(t)$ of rational functions in one variable over a field $\kk_0$, where $\kk_0$ can be any field). We will also write~$\qbar$ for the algebraic closure of~$\que$. In the number field case, we write $d = [K:\que]$ for the global degree of $K$ over $\que$; in the function field case, the global degree is $d = [K:\kk]$. When $K$ is a function field, we will distinguish two cases: we say that it is of {\it finite type $q$} if its subfield of constants is a finite field $\ff_q$ for some prime power $q$, and we say that it is of {\it infinite type} if its subfield of constants is infinite.

Next we discuss absolute values on $K$. Let $M(K)$ be the set of all places of $K$ when $K$ is a number field, and the set of all places of $K$ which are trivial over the field of constants when $K$ is a function field. For each place $v \in M(K)$ we write $K_v$ for the completion of $K$ at $v$ and let $d_v$ be the local degree of $K$ at $v$, which is $[K_v:\que_v]$ in the number field case, and $[K_v:\kk_v]$ in the function field case.

If $K$ is a number field, then for each place $v \in M(K)$ we define the absolute value $|\ |_v$ to be the unique absolute value on $K_v$ that extends either the usual absolute value on $\real$ or $\cee$ if $v | \infty$, or the usual $p$-adic absolute value on $\que_p$ if $v|p$, where $p$ is a prime. We also write $\O_K$ for the ring of integers of $K$.

If $K$ is a function field, then all absolute values on $K$ are non-archimedean. For each $v \in M(K)$, let $\OO_v$ be the valuation ring of $v$ in $K_v$ and $\MM_v$ the unique maximal ideal in $\OO_v$. We choose the unique corresponding absolute value $|\ |_v$ such that:
\begin{trivlist}
\item (i) if $1/t \in \MM_v$, then $|t|_v = e$,
\item (ii) if an irreducible polynomial $p(t) \in \MM_v$, then $|p(t)|_v = e^{-\deg(p)}$.
\end{trivlist}

\noindent
In both cases, for each non-zero $a \in K$ the {\it product formula} reads
\begin{equation}
\label{product_formula}
\prod_{v \in M(K)} |a|^{d_v}_v = 1.
\end{equation} 

We extend absolute values to vectors by defining the local heights. For each $v \in M(K)$ define a local height $H_v$ on $K_v^N$ by
$$H_v(\bx) = \max_{1 \leq i \leq N} |x_i|^{d_v}_v,$$
for each $\bx \in K_v^N$. Also, for each $v | \infty$ we define another local height
$$\H_v(\bx) = \left( \sum_{i=1}^N |x_i|_v^2 \right)^{d_v/2}.$$
Then we can define two slightly different global height functions on $K^N$:
\begin{equation}
\label{global_heights}
H(\bx) = \left( \prod_{v \in M(K)} H_v(\bx) \right)^{1/d},\ \ \H(\bx) = \left( \prod_{v \nmid \infty} H_v(\bx) \times \prod_{v | \infty} \H_v(\bx) \right)^{1/d},
\end{equation}
for each $\bx \in K^N$. These height functions are {\it homogeneous}, in the sense that they are defined on the projective space $\pee^{N-1}(K)$ thanks to the product formula (\ref{product_formula}): $H(a \bx) = H(\bx)$ and $\H(a \bx) = \H(\bx)$ for any $\bx \in K^N$ and $0 \neq a \in K$. It is easy to see that
$$H(\bx) \leq \H(\bx) \leq \sqrt{N} H(\bx).$$
Notice that in case $K$ is a function field, $M(K)$ contains no archimedean places, and so $H(\bx) = \H(\bx)$ for all $\bx \in K^N$. We also define the {\it inhomogeneous} height
$$h(\bx) = H(1,\bx),$$
which generalizes Weil height on algebraic numbers: for each $\alpha \in K$, define
$$h(\alpha) = \prod_{v \in M(K)} \max \{ 1, |\alpha|_v \}^{d_v/d}.$$
Clearly, $h(\bx) \geq H(\bx)$ for each $\bx \in K^N$. All our inequalities will use heights $H$ and $h$ for vectors, however we use $\H$ to define the conventional Schmidt height on subspaces in the manner described below.
\smallskip

We extend the height $H$ to polynomials (in particular, to quadratic forms) by viewing it as height function of the coefficient vector of a given polynomial. We also define a height function on subspaces of $K^N$. Let $V $ be an $L$-dimensional subspace of $K^N$, $1 \leq L \leq N$, and choose a basis $\bx_1,\dots,\bx_L$ for $V$ over $K$. The wedge product $\bx_1 \wedge \dots \wedge \bx_L$ can be identified with the vector of Pl\"ucker coordinates of $V$, viewed under the canonical embedding into $K^{\binom{N}{L}}$. Then we can define
$$H(V) = \H(\bx_1 \wedge \dots \wedge \bx_L).$$
The product formula implies that this definition does not depend on the choice of the basis for $V$.

An important observation is that the normalizing exponent $1/d$ in (\ref{global_heights}) makes our heights {\it absolute}, meaning that they do not depend on the number field or function field of definition, and hence are well defined over the algebraic closure of~$K$. A crucial property of height functions over number fields and function fields of finite type, which makes them the ``right tool" from the stand point of search bounds is Northcott's finiteness property (see~\cite{bombieri}, pages 25, 44, 117, 298 for a detailed discussion of Northcott's theorem and Northcott's property):

{\it Let $K$ be a number field or a function field of finite type, $N$ a positive integer, and $C$ a positive real number. Then the sets
$$\{ [\bx] \in \pee^{N-1}(K) : H(\bx) \leq C\},\ \{ \bx \in K^N : h(\bx) \leq C \}$$
are finite.}

\begin{rem} \label{heights_rem} Extensive accounts of the theory of height functions in the context of Diophantine problems can be found in~\cite{bombieri}, \cite{lang_dioph}, and~\cite{hindry_silverman}.
\end{rem}
\bigskip

We also introduce some basic language of quadratic forms (see, for instance, Chapter 1 of \cite{scharlau}, as well as~\cite{cassels_book} and~\cite{omeara} for an introduction into the subject). We write
$$F(\bX,\bY) = \sum_{i=1}^N \sum_{j=1}^N f_{ij} X_i Y_j$$
for a symmetric bilinear form in $2N$ variables with coefficients $f_{ij} = f_{ji}$ in $K$, and $F(\bX) = F(\bX,\bX)$ for the associated quadratic form in $N$ variables. We write $H(F)$ for the height of $F$, which is the height of its coefficient vector, as specified above for polynomials. Let $V \subseteq K^N$ be an $L$-dimensional subspace, $2 \leq L \leq N$, then $F$ is also defined on $V$, and we write $(V,F)$ for the corresponding quadratic space. 

A point $\bx$ in a subspace $U$ of $V$ is called singular if $F(\bx,\bwy) = 0$ for all $\bwy \in U$, and it is called nonsingular otherwise. For a subspace $U$ of $(V,F)$, define its radical 
$$U^{\perp} := \{ \bx \in U : F(\bx, \bwy) = 0\ \forall\ \bwy \in U \}$$
to be the space of all singular points in $U$. We define $\lambda(U):=\dim_K U^{\perp}$, and will write $\lambda$ to denote $\lambda(V)$. A subspace $U$ of $(V,F)$ is called regular if $\lambda(U)=0$.

A point $\bo \neq \bx \in V$ is called isotropic if $F(\bx)=0$ and anisotropic otherwise. A subspace $U$ of $V$ is called isotropic if it contains an isotropic point, and it is called anisotropic otherwise. A totally isotropic subspace $W$ of $(V,F)$ is a subspace such that for all $\bx,\bwy \in W$, $F(\bx,\bwy)=0$. All maximal totally isotropic subspaces of $(V,F)$ contain $V^{\perp}$ and have the same dimension. Given any maximal totally isotropic subspace $W$ of $V$, we define the Witt index of $(V,F)$ to be 
$$\W = \W(V) := \dim_K(W)-\lambda.$$
If $K=\kbar$, then $\W= [(L-\lambda)/2]$, where $[\ ]$ stands for the integer part function.

If two subspaces $U_1$ and $U_2$ of $(V,F)$ are orthogonal, we write $U_1 \perp U_2$ for their orthogonal sum. If $U$ is a regular subspace of $(V,F)$, then $V = U \perp \left( \perp_V(U) \right)$ and $U \cap  \left( \perp_V(U) \right) = \{\boldsymbol 0\}$, where
\begin{equation}
\label{perp_V}
\perp_V(U) := \{ \bx \in V : F(\bx,\bwy) = 0\ \forall\ \bwy \in U \}
\end{equation}
is the orthogonal complement of $U$ in $V$. Two vectors $\bx,\bwy \in V$ are called a hyperbolic pair if $F(\bx) = F(\bwy) = 0$ and $F(\bx,\bwy) \neq 0$; the subspace $\hyp(\bx,\bwy) := \spn_K \{\bx,\bwy\}$ that they generate is regular and is called a hyperbolic plane. An orthogonal sum of hyperbolic planes is called a hyperbolic space. Every hyperbolic space is regular. It is well known that there exists an orthogonal Witt decomposition of the quadratic space $(V,F)$ of the form
\begin{equation}
\label{decompose}
V = V^{\perp} \perp \hyp_1 \perp\ \dots \perp \hyp_{\W} \perp U,
\end{equation}
where $\hyp_1, \dots, \hyp_{\W}$ are hyperbolic planes and $U$ is an anisotropic subspace, which is determined uniquely up to isometry. The rank of $F$ on $V$ is $r:=L-\lambda$. In case $K=\kbar$, $\dim_K U = 1$ if $r$ is odd and 0 if $r$ is even. Therefore a regular even-dimensional quadratic space over $\qbar$ is always hyperbolic.

A nonsingular linear map $\sigma : V \to V$ is called an isometry of the quadratic space $(V,F)$ if $F(\sigma(\bx),\sigma(\bwy)) = F(\bx,\bwy)$ for all $\bx,\bwy \in V$. The set of all isometries of $(V,F)$ forms a group, denoted $\O(V,F)$. An isometry of $(V,F)$ whose set of fixed points is a co-dimension one subspace of $V$ is called a reflection. The celebrated Cartan-Dieudonn\'e theorem (see, for instance~\cite{omeara}) states that every isometry of a regular $L$-dimensional quadratic space over a field of characteristic $\neq 2$ can be represented as a product of at most $L$ reflections. Finally, we also mention the classical notion of integral equivalence between different quadratic forms: two integral quadratic forms $F$ and $G$ in $N \geq 2$ variables are said to be integrally equivalent if there exists a matrix $A \in \GL_N(\zed)$ such that $F(A\bX)=G(\bX)$.

We are now ready to proceed.
\bigskip

\section{Extensions over global fields}
\label{fields}

The first general version of Cassels' theorem over number fields has been proved by S. Raghavan in 1975,~\cite{raghavan} (bounds for diagonal ternary quadratic forms over number fields have previously been given by C. L. Siegel~\cite{siegel_ternary}). Raghavan's result is stated in terms of a slightly different height function than those we introduced, which is defined by taking the maximum over the archimedean absolute values instead of a product, and features the same exponent $(N-1)/2$ as the original Cassels' bound. His proof extends the argument of Davenport~\cite{davenport1}. In the same paper, Raghavan also produced an analogous result for zeros of hermitian forms over number fields, although it is not clear whether his bound in this case is sharp.

The first instances of Cassels'-type result for quadratic forms over functions fields  were produced by A. Prestel~\cite{prestel} in 1987 for rational function fields and by A. Pfister~\cite{pfister} in 1997 for algebraic function fields. Prestel's result is for rational function fields with any coefficient field of characteristic $\neq 2$, and the height function used is an additive (logarithmic) analogue of the projective height $H$ we defined above. His argument follows Cassels' argument \cite{cassels_book}. Prestel also establishes an interesting fact that a Cassels'-type bound cannot in general exist over rational function fields in more than one variable. Pfister's result works over algebraic function fields with arbitrary coefficient fields, and is stated in terms of the height which is the degree of the pole divisor of a point (i.e., an additive version of inhomogeneous height). Pfister's argument uses the Riemann-Roch theorem to replace the Euclidean algorithm. The exponent on the height of quadratic form (if written multiplicatively) in both of these results is $(N-1)/2$, same as Cassels'.

As a direct implication of standard height inequalities, all of these results can be combined into the following unified version of Cassels' theorem with respect to the most commonly used height functions as defined above.

\begin{thm} \label{cas1} Let $K$ be a number field or a function field, and let $F(\bX)$ be an isotropic quadratic form in $N$ variables over $K$. Then there exists $\bx \in K^N$ such that $F(\bx)=0$ and
$$h(\bx) \ll_{K,N} H(F)^{\frac{N-1}{2}}.$$
\end{thm}

As can be expected, an analogous (much stronger) bound over an algebraically closed field is easy to obtain. For instance, in~\cite{me:quad} the bound $\ll_N H(F)^{1/2}$ for the height of a nontrivial zero of an isotropic quadratic form in $N \geq 2$ variables over $\qbar$ is established.
\bigskip

\section{Multiple zeros and isotropic subspaces}
\label{subspaces}

The next natural extension of Cassels' theorem is existence of a collection of multiple linearly independent small-height zeros of a given isotropic quadratic form over a fixed field. Suppose that $F$ is a rational quadratic form in $N \geq 2$ variables, which has nontrivial isotropic points in a lattice $\Lambda \subset \real^N$. Birch and Davenport~\cite{birch_davenport} extended the argument of~\cite{davenport1} to generalize Cassels' theorem in the following way: there exists a nonzero isotropic point $\bx$ of $F$ in $\Lambda$ such that
\begin{equation}
\label{b-d}
|\bx| \ll_N |F|^{\frac{N-1}{2}} \det(\Lambda),
\end{equation}
where $|\bx|,|\bwy|,|F|$ stand for the sup-norms of the vectors $\bx$, $\bwy$, and the coefficient vector of~$F$. Further, in his posthumous 1971 paper ~\cite{davenport} (prepared by D. J. Lewis) Davenport used a geometric argument to establish the existence of a linearly independent pair of zeros $\bx,\bwy$ of $F$ in $\Lambda$ so that
\begin{equation}
\label{dav}
|\bx| \cdot |\bwy| \ll_N |F|^{N-1} \det(\Lambda)^2,
\end{equation}
This result has been extended over number fields by Chalk~\cite{chalk} in 1980 (see also~\cite{chalk1} for an earlier announcement of Chalk's result). More precisely, for a quadratic form $F$ over a number field $K$ (of degree $d$) which is isotropic over some order $\O$ in $K$, Chalk established the existence of a linearly independent pair of zeros of $F$ in $\O^N$ with their heights bounded analogously to~\eqref{dav} in terms of the height of $F$ with $\det(\Lambda)^2$ replaced by $\disc(\O)^{N/d}$ (Chalk used the same height as Raghavan~\cite{raghavan}). 

More generally, by analogy with Minkowski's successive minima theorem, one may wonder if there exist linearly independent zeros $\bx_1,\dots,\bx_N \in \Lambda$ of a rational quadratic form $F$ in $N \geq 3$ variables, isotropic on the lattice $\Lambda$, satisfying
\begin{equation}
\label{s-p}
|\bx_1| \cdots |\bx_N| \ll_N |F|^{\frac{N(N-1)}{2}} \det(\Lambda)^N.
\end{equation}
In \cite{schulze-pillot}, Schulze-Pillot showed that in fact such an estimate cannot hold for a general lattice $\Lambda$, however establishes a similar inequality with the upper bound $\ll_N |F|^{\frac{N^2}{2}-1}$ in the special case $\Lambda=\zed^N$. It has then been pointed out by W. M. Schmidt that the exponent in Schulze-Pillot's bound is sharp at least for $N=3,4,5$. Moreover, Schulze-Pillot proved that for any lattice $\Lambda$ the following inequality holds for a collection of $N$ linearly independent zeros of $F$ in $\Lambda$:
\begin{equation}
\label{s-p-1}
|\bx_1|^{N-1} \cdots |\bx_N| \ll_N |F|^{(N-1)^2} \det(\Lambda)^{2(N-1)}.
\end{equation}
A stronger estimate is possible under an additional assumption on the lattice~$\Lambda$: suppose that $\Lambda$ contains $L$ linearly independent points $\bx_1,\dots,\bx_L$, $1 \leq L \leq N$, such that $F$ vanishes identically on $\spn_{\que} \{\bx_1,\dots,\bx_L\}$, the subspace of $\que^N$ spanned by these points. In \cite{schlickewei}, Schlickewei showed that in this case there exist such points satisfying
\begin{equation}
\label{schlick}
|\bx_1| \cdots |\bx_N| \ll_N |F|^{\frac{N-L}{2}} \det(\Lambda).
\end{equation}
Schlickewei's method relies on techniques from the geometry of numbers. In~\cite{schmidt1}, W. M. Schmidt showed that Schlickewei's bound is sharp in general.

Schlickewei's result has been substantially generalized and extended over an arbitrary number field by J. D. Vaaler~\cite{vaaler:smallzeros} and over a function field by H. Locher~\cite{locher}. In fact, Vaaler's theorem provides a bound on the height of a maximal totally isotropic subspace of a quadratic  space in terms of the heights of the quadratic form and the vector space; this implies Schlickewei-type theorem as a corollary by a direct application of Siegel's lemma~\cite{vaaler:siegel}. A version of Vaaler's result over a function field of finite type, following the same method, has recently been developed in~\cite{cfh}. This method, similarly to Cassels' original argument, relies on Northcott's finiteness property, namely on our ability to choose a point or subspace over $K$ of {\it minimal height}. This approach no longer works over $\qbar$, where a version of Vaaler's result (albeit with weaker bounds) has been obtained in~\cite{me:quad} by an application of arithmetic Bezout's theorem. We record here a general version of a Cassels'-type theorem for a totally isotropic subspace of a quadratic space with the use of heights as defined in Section~\ref{heights} above.

\begin{thm} \label{smallspace} Let $K$ be a number field, a function field of finite type, or $\qbar$, and let $F$ be a nonzero quadratic form in $N \geq 2$ variables over $K$. Let $V \subseteq K^N$ be an $M$-dimensional vector space, $1 \leq M \leq N$. Let  $L \geq 1$ be the dimension of a maximal totally isotropic subspace of the quadratic space $(V,F)$, and assume that $L$ is greater than $\lambda$, the dimension of the radical of $(V,F)$. Then there exists a maximal totally isotropic subspace $U \subseteq V$ such that
\begin{equation}
\label{iso_bound}
H(U) \ll_{K,M,L,\lambda} \left\{ \begin{array}{ll}
H(F)^{\frac{M - L}{2}} H(V) & \mbox{if $K \neq \qbar$} \\
H(F)^{(L-\lambda)(L-\lambda+1)} H(V)^{\frac{4(L-\lambda)}{3}+2} & \mbox{if $K=\qbar$.}
\end{array}
\right.
\end{equation}
\end{thm}

\noindent
It should be remarked that a similar bound for the height of a totally isotropic subspace of any dimension $l < L$ follows from Theorem~\ref{smallspace} by a direct application of Siegel's lemma (cf.~\cite{vaaler:siegel},~\cite{absolute:siegel}). Indeed, let $U$ be the maximal totally isotropic subspace of $(V,F)$ satisfying~\eqref{iso_bound}, and let $\bx_1,\dots,\bx_L$ (written in the order of increasing height) be a basis for $U$ satisfying Siegel's lemma, i.e.
$$h(\bx_1) \cdots h(\bx_L) \ll_{K,L} H(U).$$
Then $U_l = \spn_K \{\bx_1,\dots,\bx_l\}$ is a totally isotropic subspace of $(V,F)$ of dimension~$l$, and
$$H(U_l) \leq h(\bx_1) \cdots h(\bx_l) \ll_{K,L,l} H(U)^{l/L},$$
so a bound on the height of $U_l$ in terms of $H(F)$ and $H(V)$ follows from~\eqref{iso_bound}.

\bigskip

\section{Effective structure theorems}
\label{structure}

Further developments in the effective theory of quadratic forms via heights included results on existence of a small-height {\it spanning family} of totally isotropic subspaces of a quadratic space. The first such result was established by H. P. Schlickewei and W. M, Schmidt~\cite{schmidt:schlickewei} over~$\que$, and then generalized over number fields by J. D. Vaaler~\cite{vaaler:smallzeros2}, building on the authors' previous results which culminated in the number field version of Theorem~\ref{smallspace}. Let $(V,F)$ be an $M$-dimensional quadratic space in $N$ variables over a number field $K$, as above, let $\lambda = \dim_K (V^{\perp})$, and $L > \lambda$ dimension of a maximal totally isotropic subspace of $(V,F)$. The Schlickewei-Schmidt-Vaaler theorem then asserts that for any $l$ with $\lambda < l \leq L$, there exist $M-l+1$ distinct $l$-dimensional totally isotropic subspaces $U_0,U_1,\dots,U_{M-l}$ of $V$ such that:
\begin{enumerate}

\item $\dim_K (U_0 \cap U_j) = l-1$ for every $1 \leq j \leq M-l$

\item $V = \spn_K \left\{ U_0 \cup \dots U_{M-l} \right\}$

\item $H(U_0)^2 \leq H(U_0) H(U_j) \ll_{K,M,l} H(F)^{M-l} H(V)^2$ for every $1 \leq j \leq M-l$.

\end{enumerate}

\noindent
In the particular case, when $\lambda=0$ and $L\geq 1$, we can take $l=1$ and conclude that there must exist a basis $\bx_0,\dots,\bx_{M-1}$ for $V$ consisting of isotropic points of $F$ so that
\begin{equation}
\label{vaa_1}
H(\bx_0)H(\bx_j) \ll_{K,M} H(F)^{M-1} H(V)^2,
\end{equation}
for every $1 \leq j \leq M-1$, and
\begin{equation}
\label{vaa_2}
H(\bx_0)^{M-1} H(\bx_1) \cdots H(\bx_{M-1})  \ll_{K,M} H(F)^{(M-1)^2} H(V)^{2(M-1)}.
\end{equation}
These bounds present generalizations of the Davenport-Chalk and Schulze-Pillot's results~\eqref{dav} and~\eqref{s-p-1}, respectively.

The first non-commutative version of a Cassels'-type result was obtained in~\cite{quaternion}. Specifically, let $D$ be a positive definite quaternion algebra over a totally real number field $K$, $F(\bX,\bY)$ a hermitian form in $2N$ variables over $D$, and $V$ a right $D$-vector space which is isotropic with respect to $F$. It is then established in~\cite{quaternion} that there exists a small-height basis for $V$ over $D$, such that $F(\bX,\bX)$ vanishes at each of the basis vectors. This result is a generalization of Vaaler's bounds~\eqref{vaa_1} and~\eqref{vaa_2} over quaternion algebras. The height functions used here in the context of quaternion algebras were first introduced by C. Lienbend\"orfer~\cite{liebendorf:1} in the case $K=\que$, and then extended in~\cite{quaternion} to any totally real number field~$K$. The main tool developed in~\cite{quaternion} is a collection of height comparison lemmas between heights over $K$ and heights over $D$. These lemmas are then applied to Vaaler's results to ``transfer" them to the quaternion algebra setting via an appropriate $K$-vector space isomorphism. In fact, this technique has further applications, as we mention below.

Theorem~\ref{smallspace} above has further applications as well, to an effective version of Witt decomposition as in~\eqref{decompose}. The following result was established over number fields in~\cite{witt}, over function fields of finite type in~\cite{cfh}, and over~$\qbar$ in~\cite{me:quad}.

\begin{thm} \label{witt_decomp} Let $K$ be a number field, a function field of finite type, or $\qbar$, and let $F$ be a nonzero quadratic form in $N \geq 2$ variables over $K$. Let $V \subseteq K^N$ be an $M$-dimensional vector space, $1 \leq M \leq N$, so that the Witt index of the quadratic space $(V,F)$ is $\W \geq 1$. Let $\lambda=\dim_K (V^{\perp})$ and $r=M-\lambda$, the rank of $F$ on $V$. There exists an orthogonal decomposition of the quadratic space $(V,F)$ of the form \eqref{decompose} with all components of bounded height. Specifically,
\begin{equation}
\label{sing_height}
H(V^{\perp}) \ll_{K,M,r} \left\{ \begin{array}{ll}
H(F)^{r/2} H(V) & \mbox{if $K \neq \qbar$} \\
H(F)^r H(V)^2  & \mbox{if $K=\qbar$,}
\end{array}
\right.
\end{equation}
and
\begin{equation}
\label{height}
\max\{ H(\hyp_i),H(U) \} \ll \left\{ \begin{array}{ll}
\left( H(F)^{\frac{M+2\W}{4}} H(V) \right)^{\frac{(\W+1)(\W+2)}{2}} & \mbox{if $K \neq \qbar$} \\
\left( H(F)^{\W^2+1} H(V)^{\frac{6\W+5}{4\W+2}} \right)^{\frac{(\W+1)(\W+2)}{2} \left(\frac{3}{2}\right)^\W} & \mbox{if $K = \qbar$,}
\end{array}
\right.
\end{equation}
for all $1 \leq i \leq \W$, where the constant in the upper bound depends on $K,N,M,\W$.
\end{thm}

\noindent
Theorem~\ref{witt_decomp} has been applied in~\cite{cfh} to establish the existence of an infinite collection of spanning families of maximal totally isotropic subspaces of bounded height as discussed above, although the bounds are weaker than those of Schlickewei-Schmidt-Vaaler.

\begin{rem} \label{alt} It is also interesting to point out that more combinatorial techniques along with Siegel's lemma have been used in \cite{chile} to establish analogues of Theorems~\ref{smallspace}  and \ref{witt_decomp} for symplectic spaces over number fields, functions fields, and algebraic closures of one or the other, all at once.
\end{rem}

\begin{rem} \label{orthogonal_basis} Another kind of orthogonal decomposition for a quadratic space $(V,F)$ is given by an orthogonal basis for $V$ with respect to $F$. The existence of such a basis of bounded height over number field, function field, and $\qbar$ is proved in~\cite{witt},~\cite{cfh}, and~\cite{me:quad}, and the analogous statement for a symplectic space is obtained in~\cite{chile}. Such results can be viewed as orthogonal versions of Siegel's lemma.
\end{rem}

Another outgrowth of Cassels'-type effective results with respect to height is an effective (weak) version of Cartan-Dieudonn\'e theorem on decomposition of isometries of a quadratic space into a product of reflections. The following result was established over number fields in~\cite{witt}, over function fields of finite type in~\cite{cfh}, and over~$\qbar$ in~\cite{me:quad}.

\begin{thm} \label{CD} Let $K$ be a number field, a function field over a perfect constant field of characteristic $\neq 2$, or $\qbar$. Let $(V,F)$ be a regular quadratic space over $K$ with $V \subseteq K^N$ of dimension $M$, $1 \leq M \leq N$, $N \geq 2$. Let $\sigma$ be an element of the isometry group $\O(V,F)$. Then either $\sigma$ is the identity, or there exist an integer $1 \leq l \leq 2M-1$ and reflections $\tau_1,...,\tau_l \in \O(V,F)$ such that 
\begin{equation}
\label{CD_0}
\sigma = \tau_1 \circ \dots \circ \tau_l, 
\end{equation}
and for each $1 \leq i \leq l$,
\begin{equation}
\label{CD_1}
H(\tau_i) \ll_{K,M} \left\{ H(F)^{\frac{M}{3}} H(V)^{\frac{M}{2}} H(\sigma) \right\}^{5^{M-1}}.
\end{equation}
\end{thm}

\noindent
Height of an isometry $\sigma$ here is defined by taking an appropriate extension of $\sigma$ to an isometry of the entire space $(K^N,F)$, and then taking the height $H(A)$ of the $N \times N$ matrix $A$  of this extended isometry, viewed as a vector in $K^{N^2}$.

A related problem for polynomial bounds on integral equivalence of quadratic forms was formulated by D. W. Masser in~\cite{masser:baker}. 

\begin{conj} \label{masser} Suppose that two nonsingular integral quadratic forms $F$ and $G$ in $N \geq 3$ variables are integrally equivalent, i.e., $F(A\bX) = G(\bX)$ for some matrix $A \in \GL_N(\zed)$. Then there exists such an integral equivalence $A$ with
$$|A| \ll_N (|F|+|G|)^{p(N)},$$
for some function $p(N)$, independent of $F$ and $G$, where $|A|$ is the sup-norm of $A$ viewed as a vector in $\zed^{N^2}$.
\end{conj}

\noindent
Masser's conjecture has been proved by Dietmann in~\cite{dietmann} for all ternary forms and in~\cite{dietmann1} for a large class of forms in $N \geq 4$ variables. The techniques used blend together ideas from arithmetic geometry and analytic number theory, and are closely related to the treatment of an analogue of Cassels' theorem over the ring of integers, which we discuss next.  
\bigskip

\section{Effective results with additional conditions}
\label{inhomog}

Another incarnation of a Cassels'-type problem is finding small-height zeros of an inhomogeneous quadratic polynomial in $N \geq 2$ variables. Contrary to the homogeneous case, here it matters whether we are working over a field or a ring. A classical version of this problem was considered by C. L. Siegel in~\cite{siegel_quadratic}, where he proved that there exists an effectively computable search bound on the sup-norm of an integral zero of an inhomogeneous quadratic polynomial $Q$ with integer coefficients in terms of the sup-norm $|Q|$ of the coefficient vector. While Siegel did not explicitly compute this bound, his method leads to an exponential dependence on~$|Q|$. Better bounds were obtained by D. M. Kornhauser in 1990~\cite{kornhauser:1},~\cite{kornhauser:2} for the cases when $N=2$ and $N \geq 5$. In particular, Kornhauser showed that in the binary case polynomial bounds in general are not possible; on the other hand, he was able to obtain polynomial bounds for $N \geq 5$. Kornhauser's results were improved and extended to $N=3,4$ by R. Dietmann~\cite{dietmann}, and then further improved by T. D. Browning and R. Dietmann~\cite{browning_dietmann}. The following theorem is due to Dietmann, except for the case $N=2$, which was established by Kornhauser.

\begin{thm} \label{rainer} Consider a quadratic polynomial in $N \geq 2$ variables with integer coefficients
$$Q(\bX) = F(\bX) + L(\bX) + A,$$
where $F$ is a nonsingular integral quadratic form, $L$ is a linear form, and $A$ is an integer. Assume that $Q$ has an integral zero $\bx$. Then there exists such a zero satisfying
\begin{equation}
\label{rainer_bound}
|\bx| \ll_N \left\{ \begin{array}{ll}
(28|Q|)^{10|Q|} & \mbox{if $N=2$,} \\
|Q|^{2100} & \mbox{if $N=3$,} \\
|Q|^{84} & \mbox{if $N=4$,} \\
|Q|^{5N+19+\frac{74}{N-4}} & \mbox{if $N \geq 5$.} \\
\end{array}
\right.
\end{equation}
\end{thm}

\noindent
Dietmann relies on geometry of numbers techniques when $N=3,4$ and uses the circle method for $N \geq 5$ (see \cite{masser:baker} for a nice expository account of Dietmann's work). On the other hand, Kornhauser's method is of more elementary nature.

A more algorithmic approach to finding integral solutions of quadratic equation was presented in the paper~\cite{grunewald_segal}, entitled {\it ``How to solve a quadratic equation in integers"}, which inspired the title of D. W. Masser's paper~\cite{masser:1} {\it ``How to solve a quadratic equation in rationals"}. Masser considers an inhomogeneous quadratic polynomial $Q$ in $N \geq 2$ variables with rational coefficients which has a rational zero, and proves that it must then have a rational zero $\bx$ whose height is $\ll_N H(Q)^{\frac{N+1}{2}}$. Moreover, he shows that this bound is best possible. Massers' idea was to homogenize $Q$ by introducing a new variable $X_{N+1}$, and then use a Cassels'-type technique to prove the existence of a small-height zero $\bx$ of the resulting quadratic form $F$ in $N+1$ variables with the additional condition that $x_{N+1} \neq 0$. Working over the field $\que$, this results in a small-height rational zero for $Q$. Then Masser's theorem can be thought of as a result on the existence of a small-height rational zero of a rational quadratic form outside of the nullspace of a linear form $X_{N+1} \neq 0$. Masser's approach has been generalized and extended in~\cite{me:smallzeros} to establish a result on the existence of a small-height zero of a quadratic form $F$ in $N \geq 2$ variables over a fixed number field $K$ outside of a union of $M \geq 1$ proper subspaces of $K^N$. In case $M>1$, the bounds of~\cite{me:smallzeros} were improved in~\cite{dietmann2}, where the following result was obtained.

\begin{thm} \label{rainer_me} Let $K$ be a number field, $F$ a quadratic form in $N \geq 2$ variables over $K$, and $V_1,\dots,V_M \subset K^N$ proper subspaces of $K^N$, $M \geq 1$. Suppose that there exists a point $\bx \in K^N \setminus \bigcup_{i=1}^M V_i$ such that $F(\bx)=0$. Then there exists such a point with
$$H(\bx) \ll_{K,N,M} H(F)^{\frac{N+1}{2}}.$$
\end{thm}

\noindent
One simple application of this result recorded in~\cite{me:smallzeros} is the observation that if $F$ has a nonsingular zero over $K$, then there exists such a zero $\bx$ with $H(\bx) \ll_{K,N} H(F)^{(N-1)/2}$; in other words, Cassels' bound holds even with the additional assumption of non-singularity.

More recently, the result of Theorem~\ref{rainer_me} has been generalized to a statement about existence of zeros of a quadratic form outside of a union of varieties. Let $K$ be a number field, function field over a finite field $\ff_q$ for some odd prime power $q$, or $\qbar$, and let $N \geq 2$, $J \geq 1$ be integers. For each $1 \leq i \leq J$, let $k_i \geq 1$ be an integer and let
$$P_{i1}(X_1,\dots,X_N),\dots,P_{ik_i}(X_1,\dots,X_N)$$
be polynomials of respective degrees $m_{i1},\dots,m_{ik_i} \geq 1$. Let 
$$Z_K(P_{i1},\dots,P_{ik_i}) = \{ \bx \in K^N : P_{i1}(\bx) = \dots = P_{ik_i}(\bx) = 0 \},$$
and define 
\begin{equation}
\label{Z_K}
\Z_K = \bigcup_{i=1}^J Z_K(P_{i1},\dots,P_{ik_i}).
\end{equation}
For each $1 \leq i \leq J$ let $M_i = \max_{1 \leq j \leq k_i} m_{ij}$, and define
\begin{equation}
\label{def_M}
M = M(\Z_K) := \sum_{i=1}^J M_i.
\end{equation}
The following theorem has been proved in~\cite{cfh}.

\begin{thm} \label{miss_hyper} Let $V \subseteq K^N$ be an $L$-dimensional vector space, $1 \leq L \leq N$. Let $F$ be a quadratic form in $N$ variables defined over $K$. Let $\W$ be the Witt index of the quadratic space $(V,F)$, $\lambda$ the dimension of its radical $V^{\perp}$, $r=L-\lambda$ the rank of $F$ on $V$, and let $m = \W+\lambda$ be the dimension of a maximal totally isotropic subspace of $(V,F)$. Let
$$Z(V,F) = \left\{ \bz \in V \setminus \{\bo\} : F(\bz) = 0 \right\}.$$
Let $\Z_K$ and $M=M(\Z_K)$ be as in \eqref{Z_K}, \eqref{def_M} above. Suppose that $Z(V,F) \nsubseteq \Z_K$. Then there exist $m$ linearly independent vectors $\bx_1,\dots,\bx_m$ in $V$ over $K$ such that $\bx_1,\dots,\bx_m \in Z(V,F) \setminus \Z_K$,
\begin{equation}
\label{h_order_1}
h(\bx_1) \leq h(\bx_2) \leq \dots \leq h(\bx_m),
\end{equation}
and for each $1 \leq n \leq m$,
\begin{equation}
\label{miss_hyper_bnd}
h(\bx_n) \ll_{K,L,M} \left\{ \begin{array}{ll}
H(F)^{\frac{9L+11}{2}} H(V)^{9L+12} & \mbox{if $K \neq \qbar$} \\
H(F)^{\max \{r,29/2\}} H(V)^{30}  & \mbox{if $K=\qbar$.}
\end{array}
\right.
\end{equation}
\end{thm}

\noindent
The method of proof employs a certain specialization argument, which works over different fields at once: first the result is established for a quadratic form containing a monomial of the form $f_{ij}X_iX_j$, and then a generic quadratic form is put into such a special form by splitting off a hyperbolic plane. At all the steps of the construction the height needs to be carefully controlled, which is done by means of a variety of previous results about points and subspaces of bounded height in a quadratic space along with Siegel's lemma. A corollary of Theorem~\ref{miss_hyper} (also obtained in~\cite{cfh}) is a statement about the existence of a flag of small-height totally isotropic subspaces of the quadratic space $(V,F)$ outside of the union of varieties~$\Z_K$. Finally, partial analogues of Theorem~\ref{rainer_me} and~\ref{miss_hyper} over a positive definite quaternion algebra over a totally real number field are obtained in~\cite{me_glenn} by ``transferring" the number field results with the use of the height comparison lemmas of~\cite{quaternion}, as discussed in Section~\ref{structure} above.
\bigskip

\section{Open problems}
\label{problems}

There are a number of further directions of investigation in the general subject of interplay of heights and quadratic forms, which started with the theorem of Cassels. We mention several general open problems, which have potential of growing into research directions in their own right.

\begin{prob} \label{two_forms} Investigate an analogue of Cassels' theorem for nontrivial simultaneous zeros of two quadratic forms $F$ and $G$ in $N \geq 2$ variables over a global field $K$. In other words, is there a search bound $C_{K,N}(F,G)$ such that whenever $F$ and $G$ have a nontrivial simultaneous zero $\bx$ over $K$, there must exist such a zero with $H(\bx) \leq C_{K,N}(F,G)$?
\end{prob}

\noindent
In view of Matiyasevich's negative answer to Hilbert's tenth problem, it is unlikely that search bounds exist for an arbitrary system of quadratic forms, since any polynomial system can be reduced to a system of quadratic equations. On the other hand, there is some hope that search bounds may still be possible for a pair of quadratic forms. In fact, R. Dietmann has recently informed me of some possible progress in this direction (joint with M. Harvey) provided that the number of variables $N$ is sufficiently big.
\smallskip

\begin{prob} \label{hermitian} Analogously to the developments for quadratic and symplectic spaces described above, produce effective results with respect to height for general hermitian forms over global fields, including bounds on totally isotropic subspaces. 
\end{prob}

\noindent
The case of hermitian forms has so far received little attention. In fact, the only result in this direction I am aware of is Raghavan's theorem in~\cite{raghavan}, which is a direct analogue of Cassels' original theorem over number fields. There appear to be many further questions for hermitian forms that should be studied.
\smallskip

\begin{prob} \label{quater} Continue the investigation of hermitian (and skew-hermitian) spaces over quaternion algebras, and possibly in more general non-commutative situations. In particular, investigate bounds on totally isotropic subspaces.
\end{prob}

\noindent
The method of \cite{quaternion} only applies to obtaining bounds on the height of isolated zeros of hermitian forms over quaternion algebras. A different technique is needed to produce results about isotropic subspaces. It is also interesting to understand if analogous results can be obtained in more general situations of central simple algebras, for instance using the height machinery as developed by Watanabe~\cite{watanabe}.
\smallskip

\begin{prob} \label{masser_conj} Prove Masser's Conjecture~\ref{masser} on small-height integral equivalences between integral quadratic forms in the remaining cases. Investigate analogous questions over number fields and their rings of integers, and extend these questions to quadratic lattices over rings of algebraic integers.
\end{prob}

\noindent
The effective structure of isometries between quadratic spaces and lattices has not been investigated much beyond Dietmann's results on Masser's conjecture and the effective version of Cartan-Dieudonn\'e theorem, as explained above.
\bigskip

{\bf Acknowledgment.} I would like to thank the referee for the valuable suggestions and corrections which improved the quality of this paper.
\bigskip

\bibliographystyle{plain}  
\bibliography{cassels}        

\end{document}